\documentclass[10pt]{article}

\usepackage{grffile}
\usepackage{ctable}

\usepackage{amssymb}
\usepackage{amsmath}
\usepackage{mathrsfs}
\usepackage[dvips]{epsfig}
\usepackage[small]{caption}
\usepackage{graphicx}
\usepackage[all]{xy}

\DeclareFontFamily{OT1}{pzc}{}
\DeclareFontShape{OT1}{pzc}{m}{it}{<-> s * [1.10] pzcmi7t}{}
\DeclareMathAlphabet{\mathpzc}{OT1}{pzc}{m}{it}

\newtheorem{Lemma}{Lemma}[section]
\newtheorem{Theorem}{Theorem}
\newtheorem{Proposition}[Lemma]{Proposition}

\newtheorem{Remark}[Lemma]{Remark}

\newenvironment{Hypothesis}[1]%
  {\begin{trivlist}\item[]{\bf Hypothesis #1 }\em}{\end{trivlist}}

\newenvironment{Proof}[1][\unskip]%
 {\begin{trivlist} \item[]{\bf Proof #1. }}%
 {\hspace*{\fill}$\rule{.4\baselineskip}{.4\baselineskip}$\end{trivlist}}

 {\begin{trivlist}\item[]\textbf{Acknowledgments.}}{\end{trivlist}}

\makeatletter\@addtoreset{figure}{section}\makeatother

\makeatletter \@addtoreset{equation}{section} \makeatother

\newcommand{\R}{\mathbb{R}}

\newcommand{\Z}{\mathbb{Z}}

\newcommand{\Ns}{\mathrm{ker\,}}

\def\Re{\mathop{\mathrm{Re}}}

\newcommand{\rmO}{\mathrm{O}}
\newcommand{\rmo}{\mathrm{o}}
\newcommand{\rmd}{\mathrm{d}}
\newcommand{\rme}{\mathrm{e}}
\newcommand{\rmi}{\mathrm{i}}

\newcommand{\id}{\mathrm{\,id}\,}

\newcommand{\eps}{\varepsilon}
\newcommand{\Nl}{\mathcal{N}}
\newcommand{\K}{\mathcal{K}}
\newcommand{\That}{\widehat{\mathcal{T}}}
\newcommand{\B}{\mathcal{B}}
\newcommand{\G}{\mathcal{G}}
\newcommand{\cS}{\mathcal{S}}
\newcommand{\cL}{\mathcal{L}}
\newcommand{\e}{\mathpzc{e}}

\newcommand{\M}{\mathcal{M}}

\newcommand{\Rm}{\mathcal{R}}

\newcommand{\diag}{\operatorname{diag}}

\makeatletter
\newsavebox{\@brx}
\newcommand{\llangle}[1][]{\savebox{\@brx}{\(\m@th{#1\langle}\)}%
  \mathopen{\copy\@brx\kern-0.5\wd\@brx\usebox{\@brx}}}
\newcommand{\rrangle}[1][]{\savebox{\@brx}{\(\m@th{#1\rangle}\)}%
  \mathclose{\copy\@brx\kern-0.5\wd\@brx\usebox{\@brx}}}
\makeatother

\definecolor{Green}{rgb}{0.,0.4,0.}

\renewcommand{\geq}{\geqslant}

\newcommand{\Rmnum}[1]{\uppercase\expandafter{\romannumeral #1\relax}}

\def\XXint#1#2#3{{\setbox0=\hbox{$#1{#2#3}{\int}$}
     \vcenter{\hbox{$#2#3$}}\kern-.5\wd0}}

\newfam\bifam
\font\tenbi=cmmib10 scaled \magstep1 \font\sevenbi=cmmib10 at 11pt
\font\fivebi=cmmib10 at 6pt \textfont\bifam = \tenbi
\scriptfont\bifam = \sevenbi \scriptscriptfont\bifam= \fivebi

\ifx\pdfoutput\undefined
   \pdffalse
\else
   \pdfoutput=1
   \pdftrue
\fi
\ifpdf
   \usepackage{graphicx}
   \usepackage{epstopdf}
\else
   \usepackage{graphicx}
\fi

\begin{document}

\begin{center}

{\fontsize{17}{17}\fontfamily{cmr}\fontseries{b}\selectfont{Bifurcation to coherent structures in nonlocally coupled systems}}\\[0.2in]
Arnd Scheel and Tianyu Tao\\
\textit{\footnotesize 
University of Minnesota, School of Mathematics,   206 Church St. S.E., Minneapolis, MN 55455, USA}
\date{\small \today} 
\vspace*{0.2in}

\end{center}

\vspace*{0.2in}

\begin{abstract}
\noindent 
We show bifurcation of localized spike solutions from spatially constant states in systems of nonlocally coupled equations in the whole space. The main assumptions are a generic bifurcation of saddle-node or transcritical type for spatially constant profiles, and a symmetry and second moment condition on the convolution kernel. The results extend well known results for spots, spikes, and fronts, in locally coupled systems on the real line, and for radially symmetric profiles in higher space dimensions. Rather than relying on center manifolds, we pursue a more direct approach, deriving leading order asymptotics and Newton corrections for error terms. The key ingredient is smoothness of Fourier multipliers arising from discrepancies between nonlocal operators and their local long-wavelength approximations. 
\end{abstract}

\section{Introduction}\label{s:1}
Much of the success of modeling has been based on infinitesimal descriptions of physical laws, leading to differential and partial differential equations as models for physical processes. However, many physical interactions are inherently nonlocal, at least at a coarser modeling level, leading to nonlocal spatial coupling, as well as dependence of time evolution on history. From a dynamical systems point of view, a natural question in this context is in how far nonlocally coupled equations may behave qualitatively differently from locally coupled problems. As usual, one can approach this question from several vantage points, striving in particular to point to phenomena where nonlocality generates new phenomena, or delineate situations where nonlocal and local equations behave in analogous fashions. Our contribution here is mainly towards the latter aspect, showing that local bifurcations in nonlocal systems produce coherent structures completely analogous to local systems. We will however also comment on phenomena that are qualitatively different as a result of nonlocality. 

We will focus on a fairly simple model problem, stationary solutions to nonlocally coupled systems, solving
\begin{equation} \label{system}
U+\K\ast U = \Nl(U;\mu) ,
\end{equation}
where $U=U(x)\in\R^k$, $x\in\R^n$, $\K\ast U$ stands for matrix convolution,
\[
(\K\ast U(x))_i = \sum_{j=1}^m \int_{\R^n} \K_{i,j}(x-y)U_j(y)dy, \hspace{0.2in} 1\le i\le k,
\]
and $\Nl(U;\mu)$ encodes nonlinear terms which depend on a parameter $\mu\sim 0$. 

We are interested in the existence of spatially localized solutions $U(x)\to U_\infty$, $|x|\to\infty$ in the  prototypical setting of a transcritical bifurcation of spatially constant solutions $U(x)\equiv m$. 

In the remainder of this introduction, we shall first provide some background and motivation for this kind of question, Section \ref{s:mot}, and then give precise assumptions and results in Sections \ref{s:set} and \ref{s:res}. We collect some standard notation used throughout at the end of the introduction.

\subsection{Motivation}\label{s:mot}

\paragraph{Applications with nonlocal coupling.} Nonlocal coupling has been suggested as a more appropriate modeling assumptions in fields all across the sciences. Prominent examples are nonlocal dispersal of plant seeds in ecology \cite{plantdisp}, fractional powers of the Laplacian as limits of random walks with Levy jumps \cite{levy}, models for neural fields \cite{neuralfieldrev}, nonlocal interactions in models for Bose-Einstein condensates \cite{bose}, kinetic equations for interacting particles \cite{swarmingrev}, shallow water-wave models \cite{waterwave}, or material science \cite{matsci}. In many of these models, one is interested in spatially localized or front-like solutions, stationary, periodic,  or  propagating in time,  which we will here refer to collectively as  coherent structures. Such solutions are usually found through a traveling-wave ansatz, thus eliminating or compactifying time dependence.  In some special situations, Many of the above examples can thereby be reduced to problems of the type \eqref{system}, and we will give details for some cases in Section \ref{s:app}. 

Arguably, the most powerful results for existence and stability of coherent structures rely on formulating the existence problem as an ordinary differential equation, a method sometimes referred to as ``spatial dynamics'' \cite{sandtw} --- for locally coupled equations, in essentially one-dimensional geometries such as the real line, cylinders $\R\times\Omega$, or with radial symmetry. We will briefly overview such results from our perspective here and discuss in how far they generalize to nonlocally coupled equations, before turning to our contribution in more detail in Sections \ref{s:set}--\ref{s:res}.

\paragraph{Local coupling --- results.} 
Replacing nonlocal coupling by diffusive coupling, say $D\Delta U$, with $D$ positive, diagonal, we can consider  the elliptic system 
\begin{equation}\label{e:rd}
\Delta U + \Nl(U;\mu)=0,
\end{equation}
with $x\in\R^n$. In the simplest case of $n=1$, one can study the resulting ordinary differential equation 
\begin{equation*}
U_x=V,\qquad 
V_x=-\Nl(U;\mu),
\end{equation*}
using dynamical systems methods such as center-manifold reduction and normal form transformations, thus leading to nearly universal predictions for bifurcations of coherent structures. In addition to existence results, such methods also allow one to state rather general uniqueness and non-existence results. In higher space dimensions, $n>1$, such reduction techniques are not known to be applicable, except for the context of radial symmetry, which allows one to formulate the existence problem as a dynamical system in the radial variable $r$, 
\begin{equation*}
U_r=V,\qquad
V_r=-\frac{n-1}{r}V-\Nl(U;\mu).
\end{equation*}
Slightly adapted center manifold and normal form theory again leads to near-universal classifications of small-amplitude coherent structures; see \cite{Srad} for a comprehensive exposition of these techniques and \cite{harioo} for background on center manifolds and normal forms. 

In the simplest case of a transcritical bifurcation in the nonlinearity, with suitable additional assumptions on the linear part, one finds a reduced equation on the center manifold of the form 
\begin{align*}
u_x&=v+\rmO\left(|\mu|(|u|+|v|)+u^2+v^2\right),\\
v_x&=\mu u - u^2 + \rmO\left(v^2+(|u|+|v|)(\mu^2+u^2+v^2)\right),
\end{align*}
which at leading order, after scaling, reduces to 
\begin{equation}\label{e:tc}
u_{xx}- u + u^2=0,
\end{equation}
with an explicit nontrivial localized solution $u(x)\to 0,\ |x|\to\infty$. Using reversibility $x\mapsto -x$ one then readily shows persistence of this solution at higher order. Exploiting the characterization of center manifolds, one also obtains non-existence of other localized solutions and, in fact, a complete characterization of small bounded solutions. 

The results in \cite{Srad} extend this machinery to radially symmetric solutions in higher space dimension, leading to a leading order equation of the form 
\[
u_{rr}+\frac{n-1}{r}u_r - u + u^2=0.
\]
This type of results is also available for Turing and Hopf bifurcations \cite{Srad}. It is however not immediately applicable to the type of nonlocally coupled equations described above, with recent progress that we shall describe next.

\paragraph{Nonlocal coupling --- center manifolds.} Going back to \eqref{system}, we can in general not find an obvious formulation as a dynamical system, with the notable exception of convolution kernels with a rational Fourier transform. Consider for instance $\K(x)=\frac{1}{2}\rme^{-|x|}$, $x\in\R$, with Fourier transform $\hat{\K}(\xi)=(1+\xi^2)^{-1}$, for which we can formally write \eqref{system} as $(\id -\Delta)^{-1}U+\Nl(U;\mu)=0$, which in turn is equivalent to 
\begin{align*}
-U+W+\Nl(U;\mu)&=0,\\
W-\Delta W &=U.
\end{align*}
Under suitable assumptions, equivalent to stability assumptions made in \cite{Srad}, one can solve the first equation for $U=\Psi(W;\mu)$ and insert into the second equation, thus obtaining a local equation for $W$, 
\[
\Delta W -W + \Psi(W;\mu)=0,
\]
which is amenable to the methods from \cite{Srad}. 

The restriction to kernels with rational Fourier transform is clearly restrictive, excluding for instance Gaussians, and one naturally wonders if similar results hold outside of this class. The more recent results in \cite{FScmfd} answer this question in the affirmative, for $x\in\R$ and $\K,\K'$ exponentially localized, such that $\hat{\K}$ is analytic in a strip of the complex plane $|\Re\xi|<\eta$. 

Some smoothness of $\K$ appears to be necessary, as the counter example of $\K=\frac{1}{2}(\delta_{-1}+\delta_{1})$ shows, which produces a simple iteration 
\[
\frac{1}{2}\left(U(x+1)-2U(x)+U(x-1)\right)+\Nl(U;\mu)=0,
\]
with completely uncorrelated solutions on lattices $x\in x_0+\Z$. A finite-dimensional reduction to an ordinary differential equation here does not seem possible. 

On the other hand, exponential localization appears to be necessary. Algebraically localized kernels typically yield algebraically localized profiles but solutions of reduced differential equations would typically converge exponentially. Nevertheless, our results remove the assumption of exponential localization, at the expense of lacking a general uniqueness argument. We still reduce to the simple ordinary differential equation \eqref{e:tc} or its higher-dimensional analogue, with exponentially localized solutions, and find weaker far-field decay only at higher order in the bifurcation parameter $\mu$.

\subsection{Setup --- linear nonlocal diffusive coupling and local bifurcations of spatially constant states}\label{s:set}
Within the context described in the previous section, we are now ready to state our main assumptions and results. We start with assumptions on the linear part, keeping in mind that the nonlinearity will be assumed to be of quadratic order in $U,\mu$, then turn to assumptions on the nonlinearity, before stating our main result. 

\paragraph{Linear diffusive coupling.}

Let $I_k$ denote the identity matrix of size $k$ and consider the linearized operator $I_k + \K*$. 

\begin{Hypothesis}{(L)}We assume that  $\K$ satisfies the following properties:
\begin{enumerate}
\item \emph{localization:  }$\K$ has finite moments of order 2, that is, $\K(x),\, |x|^2\K(x) \in L^1(\R^n,\R^{k\times k})$;
\item \emph{symmetry: } $\K(x)=\K(\gamma x)$ for all $x\in\R^n$ and all $\gamma \in \Gamma\subset \mathbf{O}(n)$, a subgroup of the orthogonal matrices with 
\[
\mathrm{Fix}\, \Gamma=\left\{x\mid\gamma x=x, \mbox{ for all } \gamma\in\Gamma\right\}=\{0\};
\]
\item \emph{minimal nullspace: } $\Ns(I_k+\int\K)=\mathrm{span}\, {\e}$ for some  $ 0\neq {\e}\in \R^k$; we then choose ${\e}^*$ such that   $\Ns(I_k+\int\K^T)=\mathrm{span}\, {\e}^*$;
\item \emph{nondegenerate second moments:}  the matrix of projected second moments $S$ with entries 
\[
S_{ij}=\int x_ix_j \langle{\e}^*,\K(x) {\e}\rangle\rmd x,
\]
is positive definite;
\item \emph{invertibility for nonzero wavenumbers: } $I_k+ \int \rme^{\rmi \langle \xi,x\rangle}\K(x)\rmd x$ is invertible for all $\xi\neq 0$. 
\end{enumerate}
\end{Hypothesis}
The assumption on positive definiteness can be readily replaced by negative definiteness, simply multiplying the equation by $-1$. 
Note that first moments, $\int x\K(x)\rmd x$ vanish, since
\begin{equation}\label{e:1st}
\int x \K(x)\rmd x=\int x K(\gamma x)\rmd x= \gamma^{-1}\int y \K(y)\rmd y,
\end{equation}
for all $\gamma\in \Gamma$, hence $\int x \K(x)\rmd x=0$ since $\Gamma$ fixes the origin, only. In fact, this is the primary reason for us to require the symmetry mentioned here. Nonvanishing first moments correspond to an effective directional transport which would need to be compensated by a drift term $c\cdot \nabla U$, say, in order to find coherent structures. 

Typical examples of symmetry groups $\Gamma$ are $\Gamma=\mathbf{O}(n)$, $\Gamma=\{\mathrm{id},-\mathrm{id}\}$ or the group generated by reflections at hyperplanes, $x_j\mapsto -x_j$. 

\begin{Remark}[normalizing second moments]\label{r:rs} There exists a coordinate change $x=T_0y$ such that the transformed kernel $\tilde{\K}(y):=|\mathrm{det}\,T_0|\K(T_0y)$ satisfies
\[
\tilde{S}_{ij}=\int x_ix_j\langle {\e}^*,\tilde{\K}(x) {\e}\rangle\rmd x=2\delta_{ij}.
\]
Indeed, let $\lambda_i$, $i=1,\ldots,k$ be the eigenvalues of $S$ and $P_i$ be the associated spectral projections. Define
\[
T_0 = \sum_{i} \sqrt{2}(\lambda_i)^{-1/2} P_i, \hspace{0.5in} \tilde{M}(y) = \langle {\e}^*, \tilde{K}(x){\e} \rangle ,\hspace{0.5in}  {M}(y) = \langle {\e}^*, {K}(x){\e} \rangle .
\]
Then
\begin{align*}
\tilde{S}_{ij} &= \int x_ix_j \tilde{M}(x) \rmd x 
= \int (T_0^{-1}y)_i(T_0^{-1}y)_jM(y)\rmd y
= \sum_{k,l}T^{-1}_{0,ik}T^{-1}_{0,jl}S_{kl}
= T_0^{-1}S(T_0^{-1})^{T}\\
&=\left(\sum_\ell \sqrt{2}\lambda_\ell^{-1/2}P_\ell\right)\left(\sum_m \lambda_m P_m\right)\left(\sum_k \sqrt{2}\lambda_k^{-1/2}P_k\right)=2\sum P_\ell=2I_k.
\end{align*}
Note that the new kernel $\tilde{\K}$ possesses all the properties assumed for $\K$ in Hypothesis (L). In particular, $\tilde{K}$ is invariant under $\Gamma$. To see this, notice that $S\gamma=\gamma S$ for all $\gamma$, and conclude that $T_0\gamma=\gamma T_0$ since spectral projections commute with $\gamma\in\Gamma$, as well. As a consequence, $\tilde{K}$ is invariant. 
\end{Remark}

The assumptions can also be stated in terms of the associated Fourier determinant $\mathcal{D}$,
\[
\mathcal{D}(\xi):=\det(I_k+\widehat{\K}(\xi)).
\]
One readily finds that $\mathcal{D}$ is of class $\mathscr{C}^2$ by the assumption on second moments, and 
\[
\mathcal{D}(0)=0,\qquad  \mathcal{D}'(0)=0, \qquad  \mathcal{D}''(0)\neq 0.
\]
The characteristic function $\mathcal{D}$ was also used in \cite{FScmfd}, identifying zeros of $\mathcal{D}$ on the real axis with bounded solutions to the linear equation, and, more generally, multiplicities of zeros adding up to the dimension of a reduced center manifold as algebraic multiplicities of a center subspace. In the setup there, $\mathcal{D}$ was analytic, allowing readily for characterizing the multiplicity of roots. We here assume just enough regularity, $\mathscr{C}^2$ to make sense of a ``double root'' of $\mathcal{D}$.

\begin{Remark}[generalizing linear assumptions]\label{r:gl}
Most examples of nonlinear problems would involve a nontrivial pointwise linear part, say, $AU+\K*U$. One quickly sees that these and more general terms should be viewed as less smooth, namely Dirac-$\delta$ contributions to the matrix kernel. Whenever this principle part, say, the matrix $A$, is invertible, the system can be readily put into our form by applying $A^{-1}$. On the other hand, when this principal part possesses a kernel, our assumption of invertibility for nonzero wavenumbers would be violated asymptotically, for $\ell\to\infty$. Bifurcation solutions in such situations are not necessarily smooth. 
\end{Remark}

Yet a different interpretation would refer to the spectrum of the linear part $I_k + \K*$, given by the set of $\lambda$ for which $I_k + \K*-\lambda I_k$ is not invertible, or, equivalently, the  closure of the set of $\lambda$ such that $\det(I_k+\widehat{\K}(\xi)-\lambda I_k)= 0$ for some $\xi\in\R$. Clearly, $\lambda=0$ is in the spectrum, choosing $\xi=0$. Also, $\lambda=0$ is minimal in multiplicity in the sense that it is an eigenvalue only for $\xi=0$, and its geometric multiplicity at $\xi=0$ is minimal. Assuming that, in addition, $\lambda=0$ is algebraically simple, one readily finds that the continuation of $\lambda$ as an eigenvalue in $\xi$ is quadratic, $\lambda= \langle \mathcal{S}\xi,\xi\rangle+\ldots$, with definite symmetric matrix  $\mathcal{S}$. Such spectral information can in general be converted into heat decay estimates for $U_t=-U+\K*U$ \cite{heatnl}.

\paragraph{Transcritical bifurcation for spatially constant solutions.}
As mentioned, we assume the presence of a simple transcritical bifurcation in the nonlinearity $\Nl$. The assumptions that follow are generic and necessary for a non-degenerate bifurcation scenario; see for example \cite{chowhale}.  They single out relevant terms in systems of equations that lead to bifurcations as in the simple scalar example $\Nl(u;\mu)=\mu u - u^2$. 
 
\begin{Hypothesis} {(TC)} We assume that $\Nl=\Nl(U;\mu):\R^k \times \R \to \R^k $ satisfies the following conditions:
\begin{enumerate}
\item 
\emph{smoothness: }
$\mathscr{C}^{K}(\R^k \times \R; \R^k)$, $K=1+\ell+2$;
\item \emph{trivial solution: } $\Nl(0;\mu) = 0$ for all $\mu$;
\item \emph{criticality:} $D_U\Nl(0;0)=0$;
\item \emph{generic linear unfolding:  }
\begin{equation}
\alpha := \langle D_{\mu,U} \Nl(0;0){\e}, {\e}^*\rangle \neq 0; \label{muvCoe}\end{equation}
\item \emph{generic nonlinearity: } 
\begin{equation}
\beta := \frac{1}{2}\langle D_{U,U} \Nl(0;0)[{\e},{\e}], {\e}^*\rangle \neq 0.  \label{QuadCoe}
\end{equation}
\end{enumerate}
\end{Hypothesis}

Smoothness assumptions ensure that the superposition operator $U(\cdot) \mapsto \Nl(U(\cdot);\mu)$ defined by $\Nl$ is of class $\mathscr{C}^1$ as a map on $H^\ell(\R^n,\R^k)$; see for instance \cite{runst1996sobolev}.

The choice of a transcritical setup here is for convenience and other elementary bifurcations can be treated in a similar fashion. A saddle-node bifurcation, where  $\langle D_{\mu} \Nl(0;0), {\e}^*\rangle \neq 0$,  $\langle D_{U,U} \Nl(0;0)[{\e},{\e}], {\e}^*\rangle \neq 0$, can be transformed into a transcritical bifurcation after subtracting one of the branches, $\tilde{U}=U-U_-(\mu)$ and reparameterizing, say, $\mu=\tilde{\mu}^2$. Also, more general nonlocal dependence of $\Nl$ on $U$ can be allowed. In such a case, Hypothesis (TC) applies to $\Nl$ acting on spatially constant solutions $U$.

\subsection{Bifurcation of spikes --- main result}\label{s:res}

We are now ready to state our main result. As suggested above, we would like to compare our problem to the local problem 
\begin{equation}\label{e:gs}
\Delta u - \alpha\mu u + \beta u^2=0,
\end{equation}
when $U\sim u{\e}$. For $\alpha\mu>0$, this equation possesses localized ground states of the form 
\begin{equation}\label{e:gs00}
u_\mathrm{c}(x;\mu)=-\alpha\beta^{-1}\mu u_*(\sqrt{\alpha\mu} x),
\end{equation}
where $u_*(y)$ is the (positive) ground state to 
\begin{equation}\label{e:gs0}
\Delta u -  u + u^2=0;
\end{equation}
see \cite{gs} for background information on existence of such ground states and their properties. 

\begin{Theorem}[bifurcation of spikes]\label{MainRes} Fix $n<6$ and $\ell>n/2$.  Assume Hypotheses (L) and (TC), and  recall the definition of $T_0$ from Remark \ref{r:rs} and $\alpha,\beta$ from \eqref{muvCoe} and \eqref{QuadCoe}.

There then exists a constant $\mu_0>0$ such that for all $0< \alpha\mu<\mu_0$, the nonlocally coupled system \eqref{system} possesses a  family of nontrivial solutions $U_*=U_*(\cdot;\mu) \in H^\ell(\R^n;\R^k)$. Moreover, $U_*(x;\mu)$ is given to leading order through
\begin{equation}\label{expan}
U_*(x;\mu) = {-\alpha\beta^{-1}\mu\left[ u_*(\sqrt{\alpha\mu} \,T_0 x)+w(\sqrt{\alpha\mu} \, x;\mu)\right]{\e} }+ u_{\perp}(x;\mu),
\end{equation}
where 
\begin{itemize}
\item $u_*$ from \eqref{e:gs00} is the  (scalar) radially symmetric ground state to \eqref{e:gs0};
\item the corrector ${w(\sqrt{\alpha\mu} \, x;\mu)}$ satisfies  $\|w(\cdot;\mu)\|_{H^\ell} \to 0$ as $\mu \to 0$;
\item $\langle u_{\perp}(x;\mu),{\e}\rangle=0$, and $\|u_{\perp}\|_{H^\ell} = \rmO(\mu^2)$.
\end{itemize}
Moreover, $U_*(x;\mu)=U_*(\gamma x;\mu)$ for all $\gamma\in\Gamma$. 
%
%
%
\end{Theorem}

We comment briefly on the scope of this result and outline the main idea of proof. 

We believe that our assumptions are to some extent sharp. Second moments are necessary to define diffusive behavior and obtain the limiting ground state problem. One would suspect that weaker localization, $\hat{\K}\sim |\xi|^{2\nu}$ for $\xi\sim 0$ would lead to reduced problems based on the fractional Laplacian $(-\Delta)^\nu$, with somewhat analogous results. Symmetry of the kernel is necessary to some extent to prevent drift of the resulting spikes. In other words, we find at leading order a manifold of translates of a ground state as solutions. The symmetry condition guarantees that all these solutions persist to higher order as non-degenerate fixed points in the fixed point subspace of the action of $\Gamma$ on profiles. We suspect that alternate assumptions could involve a variational structure in the problem. We will comment on possible drift and the technical ramifications in the discussion section.

The proof is eventually based on a rather direct contraction mapping principle. We prepare the equation by performing linear transformations singling out the neutral direction ${\e}$, followed by scaling and solving the equation in the complement ${\e}^\perp$. Finally, an ansatz substituting $u_*$ and a suitable corrector yields a contraction mapping for the corrector in a small neighborhood of the origin. 

Key to the argument is a precise formulation of the convergence of $-I_k+\K$ to $\Delta$, which we accomplish by carefully preconditioning our system; see equation \eqref{e:precond}.

\paragraph{Outline.} We perform coordinate changes and rescalings in Section \ref{s:2}, preparing for the proof through Lyapunov-Schmidt reduction and contraction mappings in Section \ref{s:3}. Section \ref{s:app} shows some more concrete applications of our result, and we conclude with a discussion in Section \ref{s:d}.

\paragraph{Notation.}

For a vector $u=(u_1,\ldots,u_k) \in \R^k$, we write $|u|$ to denote its usual Euclidean norm $|u| = \sum_{i=1}^{k} u_i^2$. We also use the standard multi-index notation in $\R^n$, that is we have $\alpha = (\alpha_1,\cdots,\alpha_n)$ with $\alpha_i \in \{0,1,\ldots \}$ and $\alpha! = \alpha_1!\cdots\alpha_n!$, $|\alpha|=\alpha_1+\cdots+\alpha_n$. So that 
$D^\alpha u = \dfrac{\partial^{|\alpha |} u}{\partial_{x_1}^{\alpha_1}\cdots \partial_{x_n}^{\alpha_n}}$.

We shall use the standard Sobolev spaces on $\R^n$ with values in $\R^k$, which are denoted by $W^{\ell,p}(\R^n, \R^k)$ or simply $W^{\ell,p}(\R^n)$ when $k=1$ or even $W^{\ell,p}$ whenever it is convenient to do so and does not cause confusions. For $\ell \ge 0$ and $1\le p \le \infty$
\[
W^{\ell,p}(\R^n,\R^k) := \{ u \in L^p(\R^n,\R^k): \partial^\alpha u \in L^p(\R^n,\R^k), 1\le |\alpha| \le \ell \},
\]
with norm
\[
\|u\|_{W^{\ell,p}(\R^n,\R^k)}=
\begin{cases}
\left(\sum_{1\le |\alpha| \le \ell} \|\partial^\alpha u \|_{L^p(\R^n;\R^k)}\right)^{1/p}, \hspace{0.1in}1\le p<\infty \\
\max_{1\le |\alpha| \le \ell}\|\partial^\alpha u\|_{L^\infty(\R^n;\R^k)}, \hspace{0.4in} p=\infty.
\end{cases}
\]
We use $H^\ell(\R^n,\R^k)$ to denote the space $W^{\ell,2}(\R^n,\R^k)$, we will also use $\mathscr{C}_b^\ell(\R^n,\R^k)$ to denote the space of $\ell-$times bounded continuously differentiable functions for $\ell=0,1,\ldots,\infty$.
 
Finally, we use the usual Fourier transform on $\R^n$, 
\[
\widehat{f} (\xi)= \int_{\R^n} f(x)\rme^{-\rmi \langle \xi,x\rangle}\rmd x {,}
\]
 for a Schwartz function $f$, which extends by isometry to all $f \in L^2(\R^n,\R^k)$.

\section{Normal forms and scalings}\label{s:2}
We change variables in $\R^k$ such that the linear part takes a simple form, Section \ref{s:2.1}, and introduce the long-wavelength scaling that exhibits the leading order asymptotics in Section \ref{s:2.2}.

\subsection{Normal forms on the linear part}\label{s:2.1}
We work in coordinates $y=T_0^{-1}x$, and drop tildes for the transformed kernel. Recall that first moments vanish \eqref{e:1st} and recall the definition of the operator $\mathcal{T}:=I_k+\K*$ with symbol $\That (\xi) = I_k + \widehat{{\K}} (\xi)$.  The next lemma isolates the center part of our equation in the first coordinate and factors a long-wavelength contribution of the convolution.

\begin{Lemma}\label{Lem1} There exist invertible $k \times k$ matrices $P, Q$, and a pseudo-differential operator $L$ with symbols $\widehat{L}(\xi),\widehat{L}^{-1}(\xi) \in L^\infty(\R^n, \R^{k\times k}) $ such that 
\[
\widehat{L}(\xi)P\That(\xi)Q = \diag\left\{\dfrac{|\xi|^2}{1+|\xi|^2},I_{k-1} \right\}.
\]
Moreover, in the canonical basis of $\R^k$, we can write $L$ and its symbol in matrix form
\[
L  = \begin{pmatrix}
L_\mathrm{cc} & L_\mathrm{ch} \\
L_\mathrm{hc} & L_\mathrm{hh} 
\end{pmatrix}, \hspace{0.2in}
\widehat{L}(\xi)  = \begin{pmatrix}
\widehat{L}_\mathrm{cc}(\xi) & \widehat{L}_\mathrm{ch}(\xi) \\
\widehat{L}_\mathrm{hc}(\xi) & \widehat{L}_\mathrm{hh}(\xi) 
\end{pmatrix},
\]
where terms with subscript $\mathrm{cc}$ denote a scalar, terms with subscript $\mathrm{ch}$ a $(k-1)$ dimensional row vector, terms with subscript $\mathrm{hc}$ a $(k-1)$ dimensional column vector, and terms with subscript $\mathrm{hh}$ a $(k-1)\times (k-1)$ matrix. We then have that $\widehat{L}_\mathrm{cc},\widehat{L}_\mathrm{ch},\widehat{L}_\mathrm{hc}$ are continuous and bounded in $\xi$, and $\widehat{L}_\mathrm{hc} \in L^\infty(\R^n, \R^k)$.
\end{Lemma}
\begin{Proof}
We divide our construction in two steps. We first normalize the constant-coefficient problem, $\xi=0$ in Fourier space, and expand. We then factor the leading-order Fourier multiplier.

\paragraph{Step 1.}

Since the rank of $\widehat{\mathcal{T}}(0)$ is equal to $k-1$, there exist invertible matrices $P$ and $Q$ such that
\[
P\That(0)Q = \diag\{0,I_{k-1}\}.
\]
By finiteness of second moment, the entries of $\That(\xi)$ are $\mathscr{C}^2$ functions in $\xi$. Taylor expanding near $\xi = 0$, we find 
\[
P\That(\xi)Q = \begin{pmatrix}
\That_\mathrm{cc}(\xi)& \That_\mathrm{ch}(\xi)\\
\That_\mathrm{hc}(\xi)& \That_\mathrm{hh}(\xi) 
\end{pmatrix},
\]
where $\That_\mathrm{cc} = |\xi|^2+\rmo(|\xi|^2)$ using the normalization condition $\int x_ix_j\langle {\e}^*, \tilde{\K}(x){\e}\rangle dx = 2\delta_{ij}$. The expansion for the other matrix elements follows similarly.

\paragraph{Step 2.} Set $
H(\xi) := \diag\{\frac{1+|\xi|^2}{|\xi|^2}, I_{k-1}\}$ for $\xi \neq 0$. 
We find
\begin{align*}
P\That(\xi)QH(\xi) &= \begin{pmatrix}
\That_\mathrm{cc}(\xi)& \That_\mathrm{ch}(\xi)\\
\That_\mathrm{hc}(\xi)& \That_\mathrm{hh}(\xi) 
\end{pmatrix} \begin{pmatrix}
\frac{1+|\xi|^2}{{|\xi|^2}}&0 \\
0& I_{k-1} 
\end{pmatrix}
= \begin{pmatrix}
\That_\mathrm{cc}(\xi)\frac{1+|\xi|^2}{|\xi|^2} & \That_\mathrm{ch}(\xi)\\
\That_\mathrm{hc}(\xi)\frac{1+|\xi|^2}{|\xi|^2} & \That_\mathrm{hh}(\xi)
\end{pmatrix}.
\end{align*}

The fact that $\That_\mathrm{cc}(0) = D_{\xi}\That_\mathrm{cc}(0) = 0$ and the normalization assumption on the second moment matrix together imply that $\That_\mathrm{cc}(\xi)=|\xi|^2 T_2(\xi)$ for some continuous function $T_2(\xi)$. On the other hand, we have $|\That_\mathrm{hc}(\xi)|/|\xi|^2 \le C$ for some constant $C$ near $\xi = 0$. Therefore, for $|\xi|\neq 0$ small, 
\[
P\That(\xi)QH(\xi)=
\begin{pmatrix}
1 & 0\\ 
\rmO(1) & I_{k-1}
\end{pmatrix} + \rmo(|\xi|). 
\]
It follows that $P\That(\xi)QH(\xi)$ is invertible with uniform bounds on the inverse near $\xi =0$, and its inverse is also of the form $
\begin{pmatrix}
1 & 0\\ 
\rmO(1) & I_{k-1}
\end{pmatrix} + \rmo(|\xi|)
$, for $\xi \neq 0$ small.

For $\xi$ bounded away from the origin, it follows from Hypothesis (L)(v) that $P\That(\xi)QH(\xi)$ is invertible for each $\xi$. Moreover, $\That(\xi) \to I_k$ for $|\xi|\to \infty$  and $H(\xi) \to I_k$ by Riemann-Lebesgue. We therefore conclude that $P\That(\xi)QH(\xi)$ is  invertible on $\R^n\backslash\{0\}$ with uniform bounds.

 We then define the pseudo-differential operator $L$ by setting its symbol $\widehat{L}(\xi)$ equal to $[P\That(\xi)QH(\xi)]^{-1}$. Then $\widehat{L} (\xi) \in L^\infty$ and it follows that
 \[
 L(\xi)P\That(\xi)Q = H(\xi)^{-1}\left[P\That(\xi)Q\right]^{-1}P\That(\xi)Q = H(\xi)^{-1}=\diag\left\{\dfrac{|\xi|^2}{1+|\xi|^2},I_{k-1}\right\},
 \]
 as claimed.
 
 If we write
\[
L  = \begin{pmatrix}
L_\mathrm{cc} & L_\mathrm{ch} \\
L_\mathrm{hc} & L_\mathrm{hh} 
\end{pmatrix}, \hspace{0.2in}
\widehat{L}(\xi)  = \begin{pmatrix}
\widehat{L}_\mathrm{cc}(\xi) & \widehat{L}_\mathrm{ch}(\xi) \\
\widehat{L}_\mathrm{hc}(\xi) & \widehat{L}_\mathrm{hh}(\xi) 
\end{pmatrix},
\]
then it follows from the above computations that the entries $\widehat{L}_\mathrm{cc}, \widehat{L}_\mathrm{ch},\widehat{L}_\mathrm{hh}$ are continuous at $\xi = 0$ and $\mathscr{C}^2$ for $\xi \neq 0$, while $\widehat{L}_\mathrm{hc}$ is only bounded in $\xi$. This verifies the last claim in the lemma and concludes the proof. 
 \end{Proof}

 Since $\widehat{L} \in L^\infty$, we know that the pseudo-differential operator $L$ maps $H^\ell$ into $H^\ell$ and is bounded. Denote by $M$ the pseudo-differential operator with symbol $|\xi|^2/(1+|\xi|^2)=: m(\xi)$, and set $V(y)=Q^{-1}U(y)$, with standard coordinates $V(y)=(v_\mathrm{c}(y),v_\mathrm{h}(y))^T$. Then, after precondition \eqref{system} with $LP$, we obtain the equivalent equation
\begin{equation}\label{TranEq}
{\diag{\{M, I_{n-1}\}}} V(y)=LP\Nl(QV(y);\mu).
\end{equation}
In particular, the linear part of the system is block-diagonal, decoupled between invertible and diffusive components. We will  next introduce rescalings to further simplify this equation.

\subsection{Taylor expansion and rescaling}\label{s:2.2}
We start by expanding our bifurcation equation and then introduce appropriate scalings. 
\paragraph{Taylor jets of the nonlinearity.} 
We write $v_\mathrm{h} = (v_2,\ldots,v_{k})^T$ in the standard coordinates of $\R^{k-1}$ and set $\mathcal{H}(V;\mu) :=P\Nl(QV;\mu)$. Then, with respect to the standard basis in $\R^k$, we denote by $\mathcal{H}_\mathrm{c}(V;\mu)$ the first component of the nonlinearity $\mathcal{H}$, and by $\mathcal{H}_\mathrm{h}(V;\mu)$ the remaining $k-1$ components. 

In this notation, \eqref{system} becomes
\begin{align}
M v_\mathrm{c} + L_\mathrm{cc}\mathcal{H}_\mathrm{c}(v_\mathrm{c},v_\mathrm{h};\mu) + L_\mathrm{ch}\mathcal{H}_\mathrm{h}(v_\mathrm{c},v_\mathrm{h};\mu)= 0\label{exeqnu0},\\
v_\mathrm{h} +  L_\mathrm{hc}\mathcal{H}_\mathrm{c}(v_\mathrm{c},v_\mathrm{h};\mu) + L_\mathrm{hh}\mathcal{H}_\mathrm{h}(v_\mathrm{c},v_\mathrm{h};\mu) = 0 \label{exeqnuh}.
\end{align}

By Hypothesis (TC)(i), we may use Taylor's theorem to write $\mathcal{H}_j$ as
\begin{align*}
\mathcal{H}_j(v_\mathrm{c},v_\mathrm{h};\mu) &=\left( a^j_{101} \mu v_\mathrm{c}+a^j_{011}\mu v_\mathrm{h}+a^j_{110}v_\mathrm{c}v_\mathrm{h} + a^j_{200}v_\mathrm{c}^2+a^j_{020}[v_\mathrm{h},v_\mathrm{h}] \right)+ \Rm_j(v_\mathrm{c},v_\mathrm{h};\mu)\\
&:= \B_j(v_\mathrm{c},v_\mathrm{h};\mu)+\Rm_j(v_\mathrm{c},v_\mathrm{h};\mu),
\end{align*}
where $j=c,h$, and with multi-index notation $\omega=(l,m,n)$, $|\omega|=2$, $a^j_{\omega} = \dfrac{1}{\omega !}D^{\omega} \mathcal{H}_j(0,0;0)$. The remainder $R_j$ satisfies the pointwise estimate
\begin{equation}\label{odR}
|R_j(v_\mathrm{c},v_\mathrm{h};\mu)| = |R_j(V;\mu)| = \rmO(\mu^2|V|+\mu|V|^2+|V|^3)
\end{equation}

for $(V;\mu)$ close to zero.

We are in particular interested in the terms $\mu v_\mathrm{c}$ and $v_\mathrm{c}^2$. In \eqref{exeqnu0}, the term $\mu v_\mathrm{c}$ is preconditioned by $L_\mathrm{cc}a_{101}^c+L_\mathrm{ch}a_{101}^h$, and the coefficient of $v_\mathrm{c}^2$ is preconditioned by $L_\mathrm{cc}a_{200}^c+L_\mathrm{ch}a_{200}^h$. Using Hypothesis (TC), we claim that
\begin{equation}\label{e:alp}
\alpha = a_{101}^c =\widehat{L}_\mathrm{cc}(0)a_{101}^c+\widehat{L}_\mathrm{ch}(0)a_{101}^h , \hspace{0.1in}\text{ and }\hspace{0.1in}
\beta =a_{200}^c =\widehat{L}_\mathrm{cc}(0)a_{200}^c+\widehat{L}_\mathrm{ch}(0)a_{200}^h .
\end{equation}

Indeed, to verify the first assertion, use the definition of $L$ in Lemma \ref{Lem1}.  We find $\widehat{L}_\mathrm{cc}(0)=1$ and $\widehat{L}_\mathrm{ch}(0)=(0,\ldots,0)$, thus verifying the second equality $a_{101}^c=\widehat{L}_\mathrm{cc}(0)a_{101}^c+\widehat{L}_\mathrm{ch}a_{101}^h$. 
To verify the first equality, let $e_1$ denote the standard coordinate vector $(1,0,\ldots,0)^T \in \R^k$. Then the derivative 
$a_{101}^c=\dfrac{\partial^2}{\partial \mu \partial v_\mathrm{c}}  \mathcal{H}_\mathrm{c}(0,0;0)$ is given by
\[
\langle D_{\mu,V}\mathcal{H}(0;0)e_1,e_1\rangle=\langle D_{\mu,U} P\Nl(0;0)Q e_1,e_1\rangle = \langle PD_{\mu,U}\Nl(0;0) {\e},e_1\rangle  = \langle D_{\mu,U}\Nl(0;0){\e},{\e}^*\rangle = \alpha,
\]
which verifies the first equality $\alpha=a_{101}^c$. The computations for $\beta$ are similar.

In the next paragraph, we shall make a series of rescalings to further simplify the equation and exhibit leading order parts.

\paragraph{Rescaling.} Recall that throughout $\alpha \mu > 0$. Now set $\tilde{\mu} = \alpha \mu$ and write $\sqrt{\tilde{\mu}} =: \eps$. We then rescale the functions $v_\mathrm{c}, v_\mathrm{h}$ to $\tilde{v}_\mathrm{c},\tilde{v}_\mathrm{h}$ through 
\[
\hspace{0.1in} v_\mathrm{c}(\cdot)=\frac{-1}{\beta}\eps^2\tilde{v}_\mathrm{c}(\eps \cdot), \hspace{0.1in} v_\mathrm{h}(\cdot)=\eps^2 \tilde{v}_\mathrm{h}(\eps \cdot).
\]

We substitute these variables into \eqref{exeqnu0} and \eqref{exeqnuh}, divide the first equation by $(-1/\beta)\eps^4$, the second by $\eps^2$, and then obtain,
\begin{align}
\eps^{-2}M^\eps \tilde{v}_\mathrm{c} + \sum_{j=\mathrm{c,h}} L_\mathrm{cj}^{\eps}[\tilde{\B}_j(\tilde{v}_\mathrm{c},\tilde{v}_\mathrm{h})+\eps^{-4}\tilde{\Rm}_j(\tilde{v}_\mathrm{c},\tilde{v}_\mathrm{h};\eps)]&=0,\label{rseqnu0}\\
\tilde{v}_\mathrm{h} +\sum_{j=\mathrm{c,h}} L_\mathrm{hj}^{\eps}[\eps^2\tilde{\B}_j(\tilde{v}_\mathrm{c},\tilde{v}_\mathrm{h})+\eps^{-2}\tilde{\Rm}_j(\tilde{v}_\mathrm{c},\tilde{v}_\mathrm{h};\eps)] &= 0. \label{rseqnuh}
\end{align}
Note that both \eqref{rseqnu0} and \eqref{rseqnuh} hold pointwise in $z = \sqrt{\tilde{\mu}} y$. Since $y$ is arbitrary, they hold for all $z \in \R^n$. We will subsequently view them as functional equations in $\tilde{v}_\mathrm{c}(\cdot)$ and $\tilde{v}_\mathrm{h}(\cdot)$.

In Fourier space, we find symbols for the rescaled linear operators $M^\eps $ and $L^\eps_j$, $j = \mathrm{cc,ch,hc,hh}$, of the form  $m(\eps \xi)$ and $\widehat{L}_j(\eps \xi)$, respectively. 

The rescaled nonlinear terms $\tilde{\B_j}, \tilde{\Rm_j}$ for $j=\mathrm{c,h}$ are defined through 
\begin{align*}
\tilde{\B}_j(u,v)&=\frac{a^j_{101}}{\alpha}  u+\frac{a^j_{011}}{\alpha} v+a^j_{110}uv + \frac{a^j_{200}}{-\beta}u^2+a^j_{020}(-\beta)v^2 ,\\
\tilde{R}_j ( u,v;\eps)&=\Rm_j \left(\frac{\eps^2u}{-\beta},\eps^2v;\frac{\eps^2}{\alpha}\right).
\end{align*}

In particular, the coefficient of the term $\tilde{v}_\mathrm{c}$ now equals $ a_{101}^c/\alpha=1$, and the coefficient of $\tilde{v}_\mathrm{c}^2$ now equals $a_{200}^h/(-\beta)=\beta/(-\beta)=-1$, using \eqref{e:alp}. As a consequence, we have 
\[
\tilde{B}_\mathrm{c}(\tilde{v}_\mathrm{c},\tilde{v}_\mathrm{h}) = \tilde{v}_\mathrm{c}-\tilde{v}_\mathrm{c}^2 + \rmO\left(|\tilde{v}_\mathrm{h}|+|\tilde{v} _\mathrm{h}|^2+\tilde{v}_\mathrm{c}|\tilde{v} _\mathrm{h}|\right).
\]

From the Sobolev embedding  $H^\ell(\R^n,\R^k)\hookrightarrow L^\infty(\R^n,\R^k) $ with $\ell>n/2$  and the fact that $H^\ell(\R^n,\R^k)$ is a continuous multiplication algebra, we find
\begin{equation}\label{Blest}
\|\tilde{\B}_j(u,v)\|_{H^\ell} \le C\left(\|u\|_{H^\ell}+\|v\|_{H^\ell}+\|u\|_{H^\ell}\|v\|_{H^\ell} +\|u\|_{H^\ell}^2+\|v\|^2_{H^\ell}\right),
\end{equation}
 for any $u\in H^\ell(\R^n), v\in H^\ell(\R^n,\R^{k-1})$, with some constant $C$ and  $j=\mathrm{c,h}$. For the remainder terms $\tilde{\Rm}_j$, we have 
\begin{equation} \label{Nlest}
 \|\tilde{\Rm}_j(u,v;\eps)\|_{H^\ell} = \rmO\left(\eps^6\right), \hspace{0.8in}\|D_u \tilde{\Rm}_j(u,v;\eps)\|_{H^\ell \to H^\ell} = \rmO\left(\eps^6\right),
\end{equation}
as $\eps \to 0$ from \eqref{odR}.

 For the rescaled linear operator $L^\eps$, using the definition of $L$ in Lemma \ref{Lem1} and the Fourier transform characterization of $H^\ell$, we find, for $\eps\to 0$, 
\begin{equation}\label{Liest}
\|(L^\eps_\mathrm{cc} - 1)u \|_{H^\ell} \to 0, \hspace{0.2in}\|L^\eps_\mathrm{ch}v \|_{H^\ell} \to 0, \hspace{0.2in}\|(L^\eps_\mathrm{hh}-I_{k-1})w \|_{H^\ell} \to 0, \hspace{0.2in} \text{and }\|L^\eps_\mathrm{hc}v \|_{H^\ell} \le C,
\end{equation}
where $u \in H^\ell(\R^n)$, $v\in H^\ell(\R^n,\R^{k-1})$, $w\in H^\ell(\R^n,\R^k)$, and $C$ is a constant independent of $\eps$.

We will study the behavior of the term $\eps^{-2}M^\eps v_\mathrm{c}$ as $\eps \to 0$ in the next section. 
To further ease notation, we drop tildes, and use $v_j,\B_j,\Rm_j$ ($j=\mathrm{c,h})$ for the variables also after the rescaling.

\section{Lyapunov-Schmidt reduction, leading-order ansatz, and corrections}\label{s:3}

We first solve \eqref{rseqnuh} to obtain $v_\mathrm{h}$ as a function of $v_\mathrm{c}$ using a fixed point argument in Section \ref{s:3.1}. We then substitute this function back into  \eqref{rseqnu0} to obtain a scalar equation for $v_\mathrm{c}$ and $\eps$, which we solve again using a fixed point argument in Section \ref{s:3.2}. For this, the crucial ingredients are estimates on the operator $M^\eps$ stated in Lemma \ref{estmult}. 

\subsection{Lyapunov Schmidt reduction}\label{s:3.1}

We write the left hand side of \eqref{rseqnuh} as $\G(v_\mathrm{h}; v_\mathrm{c},\eps)$, where
\[
\G(v;u,\eps) = v+ \sum_{j=c,h} L_\mathrm{hj}^{\eps}\left( \eps^2\B_j(u,v;\eps)+\eps^{-2}\Rm_j(u,v;\eps) \right). 
\]
Using estimates \eqref{Blest} and \eqref{Nlest}, we have $\G : H^\ell(\R^n,\R^{k-1}) \times H^\ell(\R^n) \to H^\ell(\R^n,\R^{k-1})$ for each $\eps >0$, small.
Note that we are treating $v_\mathrm{c}$ as an additional (Banach space-valued) parameter. The following lemma accomplishes the key reduction step to a scalar equation.
\begin{Lemma}\label{Lemuh} Fix $r>0$ not necessarily small and let $B_r$ denote the ball centered at $0$ with radius $r$ in $H^\ell(\R^n)$. Then  there exists $\eps_0>0$, sufficiently small, and a map $\psi(u,\eps): B_r \times (0,\eps_0) \to H^\ell(\R^n,\R^{k-1})$, such that $v = \psi(u, \eps)$ solves $\G(v;u,\eps) = 0$. Moreover, the map $u \mapsto \psi(u,\eps)$ is of class $\mathscr{C}^1$ for $u\in B_r$, and we have 
\[
\|\psi(u,\eps)\|_{H^\ell} = \rmO(\eps^2), \hspace{0.4in}\|D_u\psi(u,\eps)\|_{H^\ell \to H^\ell} = \rmO(\eps^2),
\] as 
$\eps \to 0$, uniformly for $u\in B_r$. Here $D_u\psi(u,\eps)$ denotes the  derivative of $\psi$ with respect to $u$ at the point $(u,\eps)$.  \end{Lemma}
\begin{Proof}We solve $\G(v;u,\eps)=0$ using a Newton iteration scheme. For $u \in B_r$ and $\eps_0$ small, we claim the following properties for $\G$:
\begin{enumerate}
\item $\|\G(0;u,\eps)\|_{H^\ell} = \rmO(\eps^2),$ uniformly in $u\in B_r$ and $\eps < \eps_0$;
\item $\G$ is smooth in $v$, and $D_v \G(0; u, \eps):H^\ell(\R^n,\R^{k-1}) \to H^\ell(\R^n,\R^{k-1})$ is bounded invertible with uniform bounds on the inverse for ${\eps}<\eps_0$ and $u \in B_r$. 
\end{enumerate}

For $(i)$, since $L^\eps$ is uniformly bounded in $\eps$, there exist a constant $C$ such that $\| L_\mathrm{hc}^\eps\|_{H^\ell \to H^\ell}+\| L_\mathrm{hh}^\eps \|_{H^\ell \to H^\ell } \le C$. We then have
\begin{align*}
\|\G(0,u;\eps)\|_{H^\ell} \le&\  \eps^2C(\|\mathcal{B}_\mathrm{c}(u,0;\eps)\|_{H^\ell} +\|\mathcal{B}_\mathrm{h}(u,0;\eps)\|_{H^\ell})\\
&\ +\eps^{-4}C(\|\Rm_\mathrm{c}(u,0;\eps)\|_{H^\ell}+\|\Rm_\mathrm{h}(u,0;\eps)\|_{H^\ell}).
\end{align*}
Using estimates \eqref{Blest} and \eqref{Nlest}, we have $\|\G(0;u,\eps)\|_{H^\ell} \le C(r) \eps^2$ uniformly in $u\in B_r$ and $\eps$ small.

For $(ii)$, we conclude the smoothness of $\mathcal{G}$ in $v$ by the smoothness of the superposition operator and the fact that $L^\eps_j$ are bounded linear operators. We compute the Fr\'echet derivative of $\G$ and obtain
\[ 
D_v\G(v;u,\eps) w = w+ \sum_{j=\mathrm{c,h}} L_\mathrm{hj}^\eps \left(\eps^2 D_v \mathcal{B}_j (u,v;\eps)+ \eps^{-2} D_v\Rm_j(u,v;\eps) \right) w,
\] 
for $w \in H^\ell(\R^n,\R^{k-1})$. Using estimate \eqref{Nlest}, we see that $D_v\G(0; u, \eps)$ is an $\rmO(\eps^2)$ perturbation of the identity as an operator on $H^\ell(\R^n,\R^{k-1})$ uniformly for $u \in B_r$. Thus, if $\eps_0$ is small enough, then for all $\eps$ with ${\eps}<\eps_0$, we have that $D_v\G(0;u,\eps)$ is bounded invertible with uniform bounds in $\eps$.

Having established $(i)$ and $(ii)$, we fix $\delta>0$ and $u \in B_r$. Let $N_\delta$ denote the closed ball of radius $\delta$ around $0$ in $H^\ell(\R^n,\R^{k-1})$, we introduce a map $\cS(\cdot; u,\eps): H^\ell(\R^n,\R^{k-1}) \to H^\ell(\R^n,\R^{k-1})$ through
\[
\cS(v; u,\eps) = v - D_v\G(0;u, \eps)^{-1}[\G(v;u,\eps)].
\]
We then find
\[
\|\cS(0;u,\eps) \|_{H^\ell} \le \|D_v\G(0;u,\eps)^{-1}\|_{H^\ell\to H^\ell} \|\G(0;u, \eps)\|_{H^\ell} = \rmO(\eps^2).
\]

Also, $D_v\cS(0;u,\eps) = 0$ by definition, and $\cS$ is smooth in $v$ by $(ii)$. Therefore, if $\delta$ is small and $v\in N_\delta$, it then follows that $\|D_vS(v;u,\eps)\|_{H^\ell \to H^\ell} \le C\delta$ for some constant $C$ independent of $\delta$.

We now start our iteration with $v_0 = 0$, $v_{n+1} := \cS(v_n;u,\eps)$, $n\ge 0$. Suppose that by induction $v_k \in N_\delta$, for $1\le k \le n$. Then
\[
\|v_{n+1}-v_n\|_{H^\ell} \le C\delta\|v_n-v_{n-1}\|_{H^\ell},
\]
by the mean value theorem. Therefore
\[
\|v_{n+1}\|_{H^\ell} \le \frac{C}{1-C\delta}\|v_1-v_0\|_{H^\ell} = \frac{C}{1-C\delta}\|\cS(0;u,\eps)\|_{H^\ell}.
\]
This implies that for $\eps$ small and $u \in B_r$, we have $v_{n+1} \in N_\delta$, and that $\cS$ is a contraction for $\delta$ sufficiently small. As in Banach's fixed point theorem, we conclude that $v_n \to v = \psi(u,\eps)$ as $n\to \infty$, and that $v$ is a fixed point of $\cS$. Note that, from the construction, we also obtain $\|\psi(u,\eps)\|_{H^\ell} = \rmO(\eps^2)$, uniformly for $u\in B_r$. 

To show the smooth dependence of $\psi(u,\eps)$ on $u$, we note that $\G(v;u,\eps)$ is also smooth in $u$ by Hypothesis (TC). Choosing $\eps$ small, the contraction constant for $\cS$ can be chosen uniformly for $u \in B_r$. Hence by the uniform contraction principle, e.g. \cite[Thm 1.244]{chicone2006ordinary}, we conclude that $\psi$ depends smoothly on $u$ as well.

Finally, in order to show that $\|D_u\psi\|_{H^\ell \to H^\ell} = \rmO(\eps^2)$, we differentiate the equation $0 = \G(\psi(u,\eps);u,\eps)$ in $u$ for $u\in B_r$ to see that $D_u\psi$ satisfies the equation
\[
D_v\G(\psi(u,\eps); u,\eps)  D_u\psi(u,\eps) + D_u\G(\psi(u,\eps);u,\eps) = 0.
\]

Now, $D_u\G(v;u,\eps)$ is of the form
\[
D_u\G(v;u,\eps) w =  \sum_{j=\mathrm{c,h}} L_\mathrm{hj}^\eps (\eps^2 D_u \mathcal{B}_j(u,v;\eps) + \eps^{-2} D_u\Rm_j(u,v;\eps) ) w.
\]
Hence, for $u \in B_r$ and $v=\psi(u,\eps) \in N_\delta$, we have $\|D_u\G(v;u,\eps)\|_{H^\ell \to H^\ell} = \rmO(\eps^2)$, again by \eqref{Blest} and \eqref{Nlest}.

On the other hand, $D_v\G(v;u,\eps)$ is uniformly invertible in $\eps$ for $v=\psi(u,\eps)\in N_\delta$ and $u\in B_r$ as shown  above. Therefore we can write $D_u\psi(u,\eps) = -[D_v\G]^{-1} D_u\G $ and conclude that 
\[
\|D_u\psi(u,\eps)\|_{H^\ell \to H^\ell} \le C(r,\delta)\eps^2.
\] This finishes the proof.
\end{Proof}

\begin{Remark} We cannot use the standard implicit function theorem directly to solve the equation $\G(v;u,\eps) = 0$ since the dependence of the convolution operators $L_\mathrm{hc}^\eps, L_\mathrm{hh}^\eps$ on $\eps$ is not well-defined at $\eps = 0$. 
\end{Remark}

\subsection{Preconditioning the reduced  equation and existence of spikes}\label{s:3.2}

We substitute $v_\mathrm{h} = \psi(v_\mathrm{c},\eps)$ from Lemma \ref{Lemuh} into  \eqref{rseqnu0} and obtain the scalar equation,
\begin{equation} \label{1dnl}
0 = \eps^{-2}M^\eps v_\mathrm{c} + \sum_{j=\mathrm{c,h}}L_\mathrm{cj}^\eps\left[B_j(v_\mathrm{c},\psi(v_\mathrm{c},\eps))+\eps^{-4}\Rm_j(v_\mathrm{c},\psi(v_\mathrm{c},\eps);\eps)\right].
\end{equation}
The key issue now is the behavior of the linear operator $M^\eps$ as $\eps \to 0$. Recall that by construction
\[
\widehat{M^\eps v}(\xi) = m(\eps \xi)\widehat{v}(\xi) = \frac{|\eps\xi|^2}{1+|\eps\xi|^2} \widehat{v}(\xi), 
\]
for any $v\in H^\ell(\R^n,\R^k)$. We then define a new operator $\mathcal{M}^\eps$ through 
\[ 
\widehat{\mathcal{M}^\eps v}(\xi) = \frac{m(\eps\xi)}{|\eps\xi|^2}\widehat{v}(\xi)=\frac{1}{1+|\eps\xi|^2} \widehat{v}(\xi). 
\] 
Since $1/(1+|\eps \xi|^2)$ is a bounded function on $\R^n$, $\M^\eps$ maps $H^\ell(\R^n,\R^k)$ into itself. 
For $v\in H^\ell(\R^n,\R^k)$, $(\M^{\eps})^{-1}$ is defined through
\[
\widehat{(\M^{\eps})^{-1}v} (\xi) = \frac{|\eps\xi|^2}{m(\eps\xi)} \widehat{v}(\xi)= (1+|\eps \xi|^2)\widehat{v}(\xi).
\]
Moreover, we have
\begin{align*}
\|((\M^\eps)^{-1}-1)v\|_{H^{\ell-2}} &=\left\| \left(1+|\eps\xi|^2-1\right)\widehat{v}(\xi)(1+|\xi|^2)^{\frac{\ell-2}{2}}\right\|_{L^2} 
\\
& \le \sup_{\ell} \left|\frac{|\eps\xi|^2}{1+|\xi|^2}\right| \|\widehat{v}(\xi)(1+\xi^2)^{\frac{\ell}{2}} \|_{L^{2}}\\ 
&\le \eps^2 \|v\|_{H^\ell}.
\end{align*}

Therefore, considered as an operator from $H^\ell (\R^n,\R^k)$ to $H^{\ell-2}(\R^n,\R^k)$, $(\M^\eps)^{-1}$ is well-defined, and we have $\|(\M^{\eps})^{-1}v - v\|_{H^{\ell-2}} \to 0$ as $\eps \to 0$ for $v \in H^\ell (\R^n,\R^k)$. This simple observation is key to identifying the leading-order terms and  we state it as a lemma.

\begin{Lemma}\label{estmult}The pseudo-differential operator $(\M^\eps)^{-1}$ with symbol $\dfrac{|\eps\xi|^2}{m(\eps\ell)}=1+|\eps\xi|^2 $ is well defined as a map from $H^\ell  (\R^n,\R^k)$ into $H^{\ell-2} (\R^n,\R^k)$. Moreover, it converges to the identity in the operator norm, 
\begin{equation}\label{e:precond}
\|(\M^\eps)^{-1}-I\|_{H^\ell \to H^{\ell-2}} = \rmO(\eps^2).
\end{equation}
\end{Lemma}
Since we are seeking solutions that inherit the symmetry of $\K$, we shall work with the subspace of functions in $H^\ell$ that are invariant under  $\Gamma \subset \mathbf{O}(n)$,
\[
H^\ell_{\Gamma}(\R^n,\R^k):=\left\{  u \in H^{\ell}(\R^n,\R^k) \ \middle|\    u(\cdot) = u (\gamma\cdot) \text{ for any } \gamma \in \Gamma\right\}.
\]
We remark that $(\M^\eps)^{-1}$ takes  $H^\ell_{\Gamma}$ into $H^{\ell-2}_{\Gamma}$ since its symbol is radially symmetric, that is, it commutes with the full group $\mathbf{O}(n)$.

We next turn to the model equation
\begin{equation}\label{e:mod}
\Delta v - v +v^2 = 0.
\end{equation}

\begin{Lemma}\label{nondeg}
Assume $n<6$ and fix $\ell\geq 2$.
Then \eqref{e:mod} possesses a unique (up to translations) localized solution $v_*(x)$, which is smooth and radially symmetric. Moreover  the linearization $\cL = -\Delta+1-2v_*$ at $v_*$ is nondegenerate in the  sense that  $\Ns \cL \cap L^2_{\Gamma}(\R^n) = \{0\}$. In particular, $\cL$ is bounded invertible from  $H^{\ell}_{\Gamma}(\R^n)$ to $H^{\ell-2}_{\Gamma}(\R^n)$.
\end{Lemma}

\begin{Proof} 
For existence and uniqueness of the ground state $v_*$, we refer to \cite[Lem. 13.3]{gs}. Next, by \cite[Lem. 13.4]{gs}, any element $\eta(x) \in \Ns \cL$ is of the form $\langle a, \nabla v_*(x)\rangle$ for some vector $a \in \R^n$. Now suppose $\eta \in \Ns \cL \cap L^2_{\Gamma}$. Then  $\eta(x) = \eta(\gamma x)$ for all $\gamma \in \Gamma$.
On the other hand, using radial symmetry of $v_*$, we have $v_*(\gamma x) ={v_*(x)}$ for all $\gamma$.  Differentiating in $x$, we have $\nabla v_* (x) = \gamma \nabla v_* (\gamma x)$. Together, this gives
\[
\langle a,\nabla v_*(x) \rangle = \eta(x) = \eta(\gamma x) = \langle a, \nabla v_* (\gamma x)\rangle = \langle  \gamma a, \nabla v_*(x)\rangle. 
\]
As a consequence,  $\gamma a - a$ is orthogonal to $\nabla v_*(x)$ for any $x \in \R^n$. However, using again the radial symmetry of $v_*$, we have $\nabla v_*(x)  = \frac{V'(|x|)}{|x|} x$ for some scalar function $V=V(r)$ defined for $r\ge 0$, not identically constant. Since $x \in \R^n$ is arbitrary, we conclude that $\gamma a =a$, which implies $a \in \mathrm{Fix}\,\Gamma = \{0\}$. Since ${-\Delta} +1$ is bounded invertible, and the multiplication operator $2v_*$ is compact as  maps from $H^\ell$ to $H^{\ell-2}$, we conclude that $\cL$ is Fredholm of index zero. Its restriction to the invariant subspace $H^\ell_\Gamma$ is therefore also Fredholm of index 0, hence invertible, as claimed. 
\end{Proof}

We are now ready to  set up the final fixed point iteration, using a preconditioning of the scalar reduced equation \eqref{1dnl} with $(\M^\eps)^{-1}$. 
\begin{Proposition}\label{prop} Assume $n<6$ and $\ell>n/2$. There is $\eps_1>0$ sufficiently small, such that for $0<\eps <\eps_1$, there exist a family of solutions to \eqref{1dnl} of the form $v_\mathrm{c}(\cdot;\eps) = v_*(\cdot)+w(\cdot; \eps)$. Here $w=w(\cdot,\eps) \in H^{\ell}_{\Gamma}(\R^n,\R)$ is a family of correctors parameterized by $\eps $ such that $\|w(\cdot,\eps)\|_{H^\ell} \to 0$ as $\eps \to 0$.
\end{Proposition}

\begin{Proof}
We substitute the ansatz $v_\mathrm{c} = v_* + w$ into \eqref{1dnl}, where $v_*$ is as stated in Lemma \ref{nondeg} and $w \in H^\ell_{\Gamma}$. We will determine an equation for $w$ and $\eps$ and show that it can be solved using a Newton iteration scheme near $(w,\eps)=(0,0)$.
First, for the term $\eps^{-2}M^\eps v_\mathrm{c}$ with $v_\mathrm{c} \in H^\ell$, we apply Fourier transform to obtain
\[
\eps^{-2}m(\eps\xi)\widehat{v}_\mathrm{c}(\xi) = -\frac{m(\eps\xi)}{|\eps\xi|^2}(-|\xi|^2)\widehat{v}_\mathrm{c}(\xi) = -\widehat{\M^\eps \Delta v_\mathrm{c}},
\]
and equation \eqref{1dnl} becomes
\[
0 = -\M^\eps \Delta v_\mathrm{c} + \left(L_\mathrm{cc}^\eps +L_\mathrm{ch}^\eps a_{101}^h\right)v_\mathrm{c}+\left(L_\mathrm{cc}^\eps a_{110}^c+L_\mathrm{ch}^\eps a_{110}^h\right)v_\mathrm{c}^2 +{\Rm(w,\psi;\eps)},
\]
where ${\Rm(w,\psi;\eps)}$ contains all the terms of order $\eps^2$ and higher,
\[
\Rm(w,\psi;\eps) =\sum_{j=\mathrm{c,h}} L_\mathrm{cj}^\eps\left[ \frac{a_{011}^j}{\alpha}\psi+a_{110}^j (v_*+w)\psi+a_{020}^j (-\beta)[\psi,\psi]+\eps^{-4}\Rm_j(v_*+w,\psi;\eps)\right].
\]
Indeed, for $w$ with $v_*+w \in B_r$, we claim that $\Rm$ satisfies the estimate $\|\Rm\|_{H^\ell} = \rmO(\eps^2)$ {and $\|D_w\Rm\|=\rmO(\eps^2)$}.  To see this, we first apply Lemma \ref{Lemuh} with $r = 2\|v_*\|_{H^\ell}$ to obtain $\psi = \psi(v_*+w,\eps)$ which satisfies $\|\psi(v_*+w,\eps)\|_{H^\ell} = \rmO(\eps^2)$.   

The linear operators $L_\mathrm{cc}^\eps, L_\mathrm{ch}^\eps$ are uniformly bounded in $\eps$, so that we have
\[
\left\|\sum_{j=\mathrm{c,h}} L_\mathrm{cj}^\eps\left( \frac{a_{011}^j}{\alpha}\psi+a_{110}^j v_\mathrm{c}\psi+a_{020}^j (-\beta)[\psi,\psi]\right)\right\|_{H^\ell} \le C(\|\psi\|_{H^\ell}+\|\psi\|_{H^\ell}^2) = \rmO(\eps^2).
\]
On the other hand, the remainders $\Rm_\mathrm{c}$ and $\Rm_\mathrm{h}$ satisfy $\|\Rm_\mathrm{c}\|_{H^\ell}= \rmO(\eps^6)$, $\|\Rm_\mathrm{h}\|_{H^\ell} = \rmO(\eps^6)$, uniformly for $v_*$ and $w$ such that $v_* +w \in B_r$ as $\eps \to 0$ by  \eqref{Blest} and \eqref{Nlest}. Therefore, we conclude that $\|\Rm(v_\mathrm{c},\psi;\eps)\|_{H^\ell} = \rmO(\eps^2)$ for $v_\mathrm{c}=v_*+w \in B_r$. { Similarly, we find that $\|D_{w}\Rm\|=\rmO(\eps^2)$ in the operator norm on $H^\ell$.}
 
Next, we add the equation $\Delta v_*-v_*+v_*^2 =0$ to the right-hand side of \eqref{1dnl} and precondition with the operator $(\M^{\eps})^{-1}$. Set $\alpha^\eps = L_\mathrm{cc}^\eps + \frac{a_{101}^h}{\alpha} L_\mathrm{ch}^\eps $, $\beta^\eps=-L_\mathrm{cc}^\eps +\frac{a_{200}^h}{-\beta} L_\mathrm{ch}^\eps $, and we find
\begin{align}
0 =&(\M^\eps)^{-1}\left[ (1-\M^\eps)\Delta v_* -\M^\eps \Delta w+\alpha^\eps(v_*+w)-v_*+\beta^\eps(v_*+w)^2+v_*^2 + \Rm \right] \nonumber \\ 
=& [(\M^{\eps})^{-1}-1]\M^\eps \Delta v_*+(\M^{\eps})^{-1}\left[ (\alpha^\eps-1)v_*+(\beta^\eps+1)v_*^2+\Rm \right]+ \nonumber \\
&-\Delta w+(\M^{\eps})^{-1}\left[\alpha^\eps w+\beta^\eps(2v_*w+w^2)\right], \nonumber \\
=: &F_1(w;\eps)+F_2(w;\eps)=: F(w;\eps).  \label{splfynl}
\end{align}

By Lemma \ref{estmult}, we have that $F$ maps $H^\ell_{\Gamma}(\R^n)$ to $H_{\Gamma}^{\ell-2}(\R^n)$. Our goal is to set up a Newton iteration scheme to solve $ F(w,\eps) =0$ for $w$ in terms of $\eps$ as a fixed point problem.

Following the strategy of Lemma \ref{Lemuh}, we shall show:
\begin{enumerate}
\item $\|F(0,\eps)\|_{H^{\ell-2}} \to 0$ as $\eps \to 0$;
\item $F(w,\eps)$ is continuously differentiable in $w$ and $D_wF(0,\eps): H^{\ell}_{\Gamma}(\R^n) \to H^{\ell-2}_{\Gamma}(\R^n)$ is uniformly invertible in $\eps$.
\end{enumerate}
For $(i)$, we note that
\[
F(0,\eps) = F_2(0;\eps) = [(\M^{\eps})^{-1}-1]\M^\eps d \Delta v_*+(\M^\eps)^{-1}[(\alpha^\eps-1)v_*+(\beta^\eps+1)v_*^2+\Rm(v_*,\psi;\eps)].
\]
From Lemma \ref{nondeg}, $\Delta v_* \in H^\ell(\R^n)$ for all $\ell$. Since $\M^\eps$ takes $H^\ell(\R^n)$ into itself and is uniformly bounded in $\eps$, we conclude from Lemma \ref{estmult} that 
$\|[(\M^{\eps})^{-1}-1]\M^\eps  \Delta v_* \|_{H^\ell} \to 0$
as $\eps \to 0$.

Moreover, by  \eqref{Liest}, it holds that $\| \alpha^\eps v -v\|_{H^\ell} \to 0$ and $\| \beta^\eps v + v\|_{H^\ell} \to 0$ as $\eps \to 0$, for any $v \in H^\ell(\R^n)$. Moreover, the remainder $\Rm(v_*,\psi;\eps)$ satisfies $\|\Rm\|_{H^\ell} = \rmO(\eps^2)$ as shown above. Hence, we conclude that 
\[
\| F(0;\eps)\|_{H^{\ell-2}} = \| F_2(0;\eps)\|_{H^{\ell-2} }\to 0,
\]
as $\eps \to 0$, which proves $(i)$.

For $(ii)$, we first verify that $F$ is continuously differentiable in $w$ from $H^\ell (\R^n)$ to $H^{\ell-2}(\R^n)$. Indeed, take $h, w_0\in H^\ell(\R^n)$ with $w_0$ fixed. We observe that $D_wF(w_0;\eps)h:H^\ell (\R^n) \to H^{\ell-2}(\R^n)$ is given by
\[
D_wF(w_0;\eps)h = -\Delta h+(\M^\eps)^{-1}\left[(a^\eps h)+2v_*\beta^\eps h + 2w_0h)+D_w\Rm h\right],
\]
which depends continuously on $w_0$ since the superposition operator induced by $\Nl$ is of class $\mathscr{C}^1$.

Now, at $w_0 = 0$, we see that, $D_wF(0;\eps)h \to -\Delta h+h-2v_*h = \cL h$ in $H^{\ell-2}(\R^n)$ as $\eps \to 0$ for $h \in H^\ell$ because $\|D_w\Rm h\|_{H^\ell} = \rmO(\eps^2)$ as noted above. By Lemma \ref{nondeg}, the operator $\cL : H^\ell_{\Gamma}(\R^n) \to H^{\ell-2}_{\Gamma}(\R^n)$ is bounded invertible. We notice that $D_wF(0;\eps)$ respects the symmetry and is a small perturbation of $\cL$, therefore invertible with uniform bounds on the inverse for $\eps$ small enough. This shows $(ii)$.

We now set up the Newton iteration scheme, by iterating the map  $\tilde{\cS}$, defined as
\[
\tilde{\cS}(w;\eps) = w-D_wF(0;\eps)^{-1}[F(w;\eps)]
.\]
Note that $\tilde{\cS}$ respects the symmetry as well, $\tilde{\cS} : H^\ell_\Gamma(\R^n) \to H^{\ell-2}_\Gamma(\R^n)$. Therefore, we can proceed as in Lemma \ref{Lemuh} to obtain $w=w(\eps)$ which solves $F(w(\eps);\eps) = 0$ for $\eps $ small enough and satisfies $\|w(\eps)\|_{H^\ell} \to 0$ as $\eps \to 0$.
\end{Proof}

Finally, we prove Theorem \ref{MainRes}.
\begin{Proof}[of Theorem \ref{MainRes}] We revert to tildes for the rescaled variables. From Proposition \ref{prop}, we know that \eqref{1dnl} has a solution of the form $\tilde{v}_\mathrm{c}(\cdot) = v_*(\cdot)+w(\cdot;\eps)$. Together with $\tilde{v}_\mathrm{h} = \psi(\tilde{v}_\mathrm{c},\eps)$, reverting the rescaling, we obtain $v_\mathrm{c}(\cdot) = -\frac{\alpha}{\beta}\mu \tilde{v}_\mathrm{c}(\sqrt{\alpha\mu }\cdot)$ and $v_\mathrm{h}(\cdot) = \alpha\mu \tilde{v}_\mathrm{h}(\sqrt{\alpha\mu}\cdot)$ as solutions to \eqref{exeqnu0} and \eqref{exeqnuh}.

Now, recall that $V=(v_\mathrm{c},v_\mathrm{h})^T$, and that the original variable $U$ is obtained as $U= QV$, where $Q$ is defined in Lemma \ref{Lem1}. We conclude that $U(\cdot)=v_\mathrm{c}(\cdot){\e}+v_{\perp}(\cdot)$, where $v_{\perp}$ takes values in the orthogonal complement of ${\e}$. The behavior of $v_\mathrm{c},v_{\perp}$ as $\mu \to 0$ is a direct consequence of Lemma \ref{Lemuh} and Proposition \ref{prop}. Finally, we restore the original variable $x = T_0y$, thus getting the desired form of the bifurcating solution.
\end{Proof}

\section{Applications}\label{s:app}

We outline how our main result applies rather immediately to a variety of specific model equations. 

\paragraph{Neural fields.}

In the simplest setup, neural field equations involve nonlocal coupling through a sigmoidal response function $S$, with dynamics 
\[
u_t=-u+\K*S(u;\mu),
\]
where $\K$ is a not necessarily positive convolution kernel sampling input from firing neighboring neurons, and $S$ is the firing rate depending on the state $u$ of neurons, typically a sigmoidal, strictly monotone function \cite{neuralfieldrev}. The sign of $\K$ may change depending on an excitatory or inhibitory coupling. Vector-valued generalizations of the equations have been proposed to model functionally different populations of neurons. 

Looking for stationary spike-like solutions of this equation, we set $u_t=0$ and substitute $U=S(u;\mu)$, with inverse $u=\Psi(U;\mu)$, obtaining
\[
-U+\K*U+(U-\Psi(U;\mu))=0,
\]
which is of the form \eqref{system}. Assuming $\int\K=1$, saddle-node bifurcations occur when the nonlinearity $U-\Psi(U;\mu)$ has a double zero, which is equivalent to a double zero in $-u+S(u;\mu)$. Hypotheses (L) and (TC) directly translate into assumptions on $\K$ and $S$. A generic saddle-node bifurcation can be easily transformed into a transcritical bifurcation as outlined in the discussion following Hypothesis (TC). 

The spikes constructed in this fashion would be expected to be unstable, of Morse index 1, their stable manifold separating spatially uniformly quiescent and spatially uniformly excited populations of neurons. 

\paragraph{Material science.} Phase separation in multi-component alloys has been modeled by free energy functionals for concentrations of species $W(u)$, $u\in \R^k$, together with a local or nonlocal interaction term. For nonlocal interactions, together with an $H^{-1}$-gradient flow, one obtains nonlocal Cahn-Morral systems
\[
u_t=-\Delta(-u+J*u-W'(u)),
\]
see \cite{cahn-morral}. Equilibria satisfy 
\[
-u+J*u-W'(u)=\mu,
\]
with chemical potential $\mu\in\R^k$. Saddle-node bifurcations in $W'(u)+\mu=0$ now lead to bifurcation of spikes as constructed here. We also note that anisotropic versions, respecting discrete crystallographic symmetries $\Gamma$, have also been proposed, at least in a context of local coupling \cite{anisoCH}. 

\paragraph{Dispersive solitary waves.}
In a slightly different direction, nonlocal coupling can encode dispersion, such as in models for shallow water waves generalizing KdV
\[
u_t=(Mu-u^2)_x,
\]
where $Mu$ is a nonlocal pseudo-differential operator generalizing $\partial_{xx}$ in the KdV equation. Traveling waves satisfy 
\[
Mu+cu-u^2=0,
\]
which in the case of $M=I_k+\K*$ reduces to the problem studied here \cite{waterwave}. Other examples include systems of nonlinear Schr\"odinger equations with nonlocal dispersion, 
\[
\rmi v_t=-v+\K*v + N(v)\in\R^k,
\]
with possible nonlocal nonlinear dispersion $N(v)$. Here, the simplest case $v=\rme^{\rmi\mu t}u,\ u\in\R$, $N(u)=u|u|^2$, leads to 
\[
-u+\K*u-\mu u + u^3=0,
\]
which is amenable to the analysis presented here, slightly changing scalings to account for the cubic nonlinearity.

%
%
%
%

\section{Discussion}\label{s:d}

We presented a direct and simple approach to the bifurcations of localized spikes in nonlocally coupled systems. While somewhat simpler and more general than approaches based on spatial dynamics, it does not offer insight into uniqueness questions, and, arguably, relies on an a priori understanding of the resulting phenomena. Our assumptions appear to be sharp in terms of localization of the convolution kernel. Interesting questions arise when studying kernels with less smoothness; see for instance \cite{pin} for a numerical exploration of kernel regularity on phenomena.

\paragraph{Large $\mu$.}
Of course, the bifurcation theoretic approach here is limited to small values of the bifurcation parameter $\mu$. Quite different phenomena are to be expected for large values of $\mu$ as the simple scalar example
\begin{equation}\label{e:simple}
-u+\K*u + \mu u - u^2=0,
\end{equation}
shows. For $\mu=1/\varepsilon$, we can scale $\varepsilon u =v$ and obtain
\[
v-v^2+\varepsilon(-v+\K*v)=0.
\]
At $\varepsilon=0$, we have solutions $v(x)=1$ for $x\in\Omega$, $v(x)=0$ otherwise, for any measurable $\Omega$. The linearization at such solutions in $L^\infty(\R^n)$ is invertible as a multiplication operator with values $\pm 1$, and the solutions therefore can be continued in $\varepsilon$. This plethora of solution does of course not exist in the case of diffusive, local coupling; see \cite{batespleth}. The transition can in some cases be understood as a depinning transition of interfaces as analyzed in \cite{pin}. 

\paragraph{Tail expansions.}
At leading order, the solutions we find here are exponentially localized, with exponential rate of order $\sqrt\mu$. As is clear from the simple example \eqref{e:simple}, the actual solution will typically not be exponentially localized. To see this, it suffices to observe that $-u+\mu u - u^2$ is exponentially localized when $u$ is, but $\K*u$ is not, for instance when both $u$ and $\K$ are positive, and $\K$ does not decay exponentially. 

Somewhat more precisely, our analysis shows that the leading order correction to the solution of 
\[
-u+\K*u-\mu u + u^2=0,
\]
can be found by substituting $u(x)=\mu u_*(\sqrt{\mu}x)+\mu w(\sqrt\mu x)$, $u_*''-u_*+u_*^2=0$,
finding at leading order 
\[
w''- w+2 u_*w=u_*''-\frac{1}{\mu}(-u_*+\K_{\sqrt{\mu}}*u_*).
\]
The right-hand side is algebraically localized with decay behavior dominated by $\K(x)\int u_*$. Solving for $ w$ and comparing decay rates shows the same behavior for the corrector $w$. In principle, one can in this way obtain higher-order algebraic expansions for the decay of $u$, assuming an expansion for $\K$ in terms of $x^{-k}$.

\paragraph{Periodic spike patterns.} The analysis of the bifurcation towards spatially periodic patterns is much simpler, due to the fact that the convolution acts as a compact perturbation of the identity. Therefore, classic Lyapunov-Schmidt bifurcation analysis will yield bifurcation of spatially periodic patterns, after restricting to appropriate symmetry planforms, such as hexagonal or square lattices in $\R^2$. More interesting and relevant for our technical questions here is the case of large spatial period, say $L= L_0/\sqrt\mu$. Imposing such boundary conditions allows for an analysis completely analogous to our analysis here, with linear Fourier multipliers restricted to periodic boundary conditions, thus allowing for the same bounds and convergence estimates. The solutions in the rescaled system would be a periodic solution to 
\[
\Delta u -u+u^2=0, 
\]
which one can obtain, using a variety of methods, for instance bifurcation and global continuation \cite{kiel}. Analyzing the maximum $A$ of the spike as a function of $L$, one then finds a stronger dependence of $A$ on $L$ for weakly decaying kernels by inspecting the residual as described above. In particular, one expects pulse interaction in the ``weak'', well-separated regime to be algebraic rather than exponential when kernels have algebraic tails. 

%

\paragraph{Traveling waves.} In some cases, relevant coherent structures may not be stationary in time, such that we need to make assumptions on temporal dynamics. Considering for instance the neural field model
\[
u_t=-u+\K*S(u),
\]
we find the traveling-wave equation 
\[
cu_x - u + \K*S(u)=0.
\]
Applying $(1-c\partial_x)^{-1}$ gives 
\[
-u+\tilde{K}_\mathrm{c}*S(u)=0,\qquad \tilde{K}_\mathrm{c}=(1-c\partial_x)^{-1}*\K.
\]
Assuming that $\K$ is of class $W^{2,1}$, say, $\tilde{\K}_\mathrm{c}$ is differentiable in $c$ as a function in $W^{1,1}$, with expansion $\tilde{\K}_\mathrm{c}=\K+c\K'$. In the long-wavelength scaling, $\K'$ converges to $\partial_x$, such that we obtain at leading order, after scaling  the model equation 
\[
u''+\tilde{c}u'+u-u^2=0.
\]
This strategy has been carried out in the case of exponentially localized kernels in \cite[\S 4.2]{FScmfd}, and we expect that the methods here allow for an adaptation to kernels with second moments, only.

\paragraph{Stability.} Focusing on the existence problem of stationary solution, we did not discuss dynamics in most of the exposition. From that perspective, a first relevant question would be the stability of bifurcating solutions. It is worth noting that an important tool for this analysis, the Evans function \cite{evans,agj,sandtw} is not available for the nonlocal equations considered here. On the other hands, an analysis as in \cite{HS2,HS}, exploiting Evans function analysis for the leading order expansion combined with a perturbation argument as presented here should yield stability information for the type of solutions constructed here for nonlocal equations. More ambitiously, it would be interesting to construct weak interaction manifolds \cite{bz,sandtw,mz}, for spikes with algebraic tail decay constructed here.

%

\end{document}